\documentclass[onefignum,onetabnum]{siamart171218}

\usepackage{tikz}
\usepackage{url}
\usepackage{amsmath,amssymb}
\usepackage[mathscr]{euscript}
\usepackage{stmaryrd}
\usepackage{mathtools}
\usepackage{algorithmic}
\Crefname{ALC@unique}{Line}{Lines} 
\newsiamthm{prop}{Proposition}
\newsiamremark{defn}{Definition}
\newsiamthm{property}{Property}
\newsiamremark{example}{Example}

\newcommand{\R}{\mathbb{R}}
\newcommand{\C}{\mathbb{C}}

\newcommand{\g}{\mathcal{G}}

\newcommand{\xv}{\mathbf{x}}
\newcommand{\uv}{\mathbf{u}}
\newcommand{\vv}{\mathbf{v}}

\newcommand{\cv}{\mathbf{c}}
\newcommand{\dv}{\mathbf{d}}

\newcommand{\zv}{\mathbf{z}}

\newcommand{\Av}{\mathbf{A}}
\newcommand{\Bv}{\mathbf{B}}
\newcommand{\Cv}{\mathbf{C}}
\newcommand{\Dv}{\mathbf{D}}

\newcommand{\Uv}{\mathbf{U}}
\newcommand{\Vv}{\mathbf{V}}

\newcommand{\Xv}{\mathbf{X}}

\ifpdf
\hypersetup{
  pdftitle={J. Pickard}
}
\fi

\newcommand{\tA}{\textsf{A}}
\newcommand{\tB}{\textsf{B}}
\newcommand{\tC}{\textsf{C}}
\newcommand{\tD}{\textsf{D}}

\newcommand{\tR}{\textsf{R}}
\newcommand{\tS}{\textsf{S}}
\newcommand{\tT}{\textsf{T}}
\newcommand{\tX}{\textsf{X}}
\newcommand{\tY}{\textsf{Y}}

\newcommand{\wv}{\mathbf{w}}
\newcommand{\yv}{\mathbf{y}}

\newcommand{\h}{\mathcal{H}}
\newcommand{\V}{\mathcal{V}}
\newcommand{\e}{\mathcal{E}}
\newcommand{\KP}{KP }

\title{Kronecker Product of Tensors and Hypergraphs: Structure and Dynamics\thanks{Submitted to the editors August 9, 2023.
    \funding{This material is based upon work supported by the Air Force Office of Scientific Research (AFOSR) under award number FA9550-22-1-0215 (IR), NSF grant DMS-2103026 
    and AFOSR grant FA9550-23-1-0400 (AMB), NIGMS GM150581 (JP), and MATHWORKS (IR). Any opinions, finding, and conclusions or recommendations expressed in this material are those of the author(s) and do not necessarily reflect the views of the United States Air Force.}
    }
}
\author{Joshua Pickard\thanks{Department of Computational Medicine \& Bioinformatics, University of Michigan, Ann Arbor, MI 48109, USA (e-mail: \email{jpic@umich.edu}).} 
\and Can Chen\thanks{School of Data Science and Society and Department of Mathematics, University of North Carolina at Chapel Hill, Chapel Hill, NC 27599, USA (\email{canc@unc.edu}).}
\and Cooper Stansbury\thanks{Department of Computational Medicine \& Bioinformatics, University of Michigan, Ann Arbor, MI 48109, USA (\email{cstansbu@umich.edu}).}
\and Amit Surana\thanks{Raytheon Technologies Research Center, East Hartford,
    CT 06108 (\email{amit.surana@rtx.com}).}
\and Anthony Bloch\thanks{Department of Mathematics, University of
    Michigan, Ann Arbor, MI 48109, USA (\email{abloch@umich.edu}).}
\and Indika Rajapakse\thanks{Department of Computational Medicine \&
    Bioinformatics, Medical School and Department of Mathematics,
    University of Michigan, Ann Arbor, MI 48109, USA (\email{indikar@umich.edu}).}
    }

\headers{Kronecker Product of Tensors and Hypergraphs}{J. Pickard et. al.}

\begin{document}

\maketitle

\begin{abstract}
Hypergraphs and tensors extend classic graph and matrix theory to account for multiway relationships, which are ubiquitous in engineering, biological, and social systems. While the Kronecker product is a potent tool for analyzing the coupling of systems in graph or matrix contexts, its utility in studying multiway interactions, such as those represented by tensors and hypergraphs, remains elusive. In this article, we present a comprehensive exploration of algebraic, structural, and spectral properties of the tensor Kronecker product. We express Tucker and tensor train decompositions and various tensor eigenvalues in terms of the tensor Kronecker product. Additionally, we utilize the tensor Kronecker product to form Kronecker hypergraphs, a tensor-based hypergraph product, and investigate the structure and stability of polynomial dynamics on Kronecker hypergraphs. Finally, we provide numerical examples to demonstrate the utility of the tensor Kronecker product in computing Z-eigenvalues, various tensor decompositions, and determining the stability of polynomial systems.
\end{abstract}


\noindent

\section{Introduction}
Engineering, biological, and physical systems are often characterized as the composite of smaller, interconnected systems. Early contributions by Kalman \cite{kalman1963mathematical} and Gilbert \cite{gilbert1963controllability} emphasized the importance of partitioning control systems, enabling the study of subsystems to gain insights into the primary system. The Kronecker product (KP) $\otimes$ is a fundamental operation that has been used across a variety of fields to describe coupled systems. \textcolor{black}{This article aims to consolidate and advance the application of the \KP in the context of tensor and hypergraphs.} 

Seeded in classic results from matrix theory \cite{zehfuss1858uber}, the \textcolor{black}{\KP has roots} in network science and dynamical systems. In the 1950s, graph products emerged as a means of investigating the combinatorial and structural properties resulting from different methods
of combining networks \cite{frucht1949groups, sabidussi1959composition, weichsel1962kronecker, vizing1963cartesian, hammack2011handbook}. Stemming from this, various graph products have since been expressed by the \KP \cite{weichsel1962kronecker, imrich2000product, hammack2011handbook}. From these structures, methods for network compression and summarization \cite{leskovec2010kronecker} along with insights to the dynamics and control of coupled systems \cite{chapman2014controllability, chapman2014kronecker, hao2019controllability} have sprouted. This framework has cultivated additional offshoots in 
multi-agent systems \cite{bizyaeva2022nonlinear, leonard2024fast}, data science \cite{van1993approximation, tsiligkaridis2013covariance, cai2022kopa}, and many other domains that utilize matrix and graphical models \cite{van2000ubiquitous}. While the \KP of matrices and graphs offer a foundation for coupling systems, recent research 
builds beyond pairwise structures of matrices and graphs towards tensors and hypergraphs \cite{kolda2009tensor, battiston2020networks}.

Hypergraphs extend graph theory by allowing hyperedges to contain more than two vertices \cite{berge1984hypergraphs}. Pairwise interactions alone fail to unambiguously capture the structure or dynamics of many systems. For instance, social and communication structures with ``friend groups" or ``group chats" cannot be adequately described by a list of \textit{pairs} of people that interact, compelling the adoption of hypergraphs to represent \textit{multiway}, group interactions. This approach has been adopted in ecology \cite{golubski2016ecological}, genomic and polymer analysis \cite{dotson2022deciphering}, and many other domains. Although hypergraph products have been utilized since the early 1990s \cite{sonntag1990hamiltonicity, sonntag1992hamiltonicity} and remain an active area of study \cite{hellmuth2012survey, ostermeier2012cartesian, kaveh2015hypergraph,hellmuth2016fast}, previous investigations into hypergraph products have been limited to a tensor-free approach.

\begin{figure}[t]
    \centering
    \includegraphics[width=\textwidth]{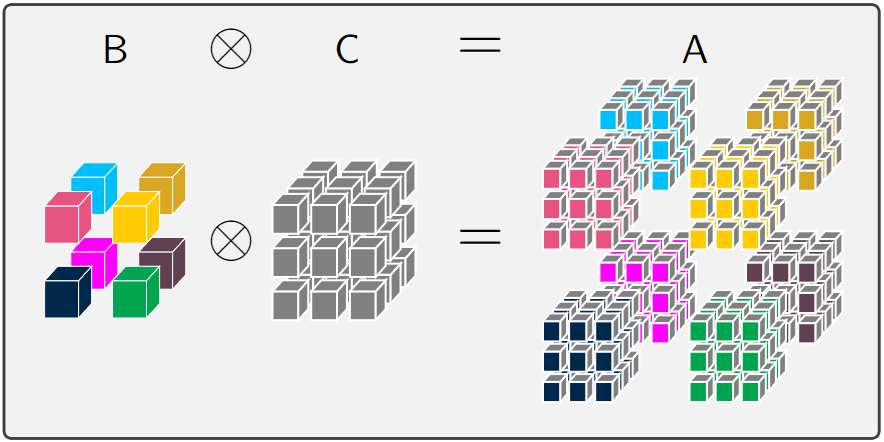}
    \caption{An example of the tensor Kronecker product.}
    \label{fig: Kronecker block tensor}
\end{figure}

Tensors are the matrices of multiway interactions, playing a crucial role in hypergraph analysis. The multiway analogue of the graph adjacency matrix is an adjacency tensor(s), offering a means to explore the structure and dynamics of the network. Tensor eigenvalues, distances, and decompositions have been utilized to study hypergraph centrality \cite{benson2019three}, similarity \cite{surana2022hypergraph}, and entropy \cite{chen2020tensor}. The homogeneous polynomials of tensors have been used to model hypergraph dynamics, from which the controllability \cite{chen2021controllability}, observability \cite{pickard2023observability, pickard2024geometric}, and stability \cite{chen2021stability, chen2022explicit, cui2024metzler} have been investigated. Tensor-based hypergraph products are imperative to quantitatively study the compositions of hypergraphs. 
Nevertheless, the tensor \KP for this purpose remains relatively unexplored.

\textcolor{black}{Tensor-based hypergraph products provide a natural means to quantitatively study hypergraph products, and several approaches to the tensor \KP (see \cref{fig: Kronecker block tensor}) have been made.}
To the best of our knowledge, the tensor \KP was first introduced by Nicholas in 2001 \cite{nicholson2001kronecker1} for 3rd and 4th order tensors in mechanics. Parallel to the matrix KP, which despite its homage to Kronecker was first published by Georg Zehfuss in 1858 \cite{zehfuss1858uber} as documented in \cite{muir1906theory, henderson1983history},  the tensor \KP operation has been reintroduced and misattributed several times. For instance, an identical operation dubbed ``recursive tensor multiplication" was defined in \cite{akoglu2008rtm}, where it is utilized to study network dynamics, and has been given attribution by \cite{eikmeier2018hyperkron, colley2023dominant}; the operation was defined again in \cite{phan2012revealing} for identifying patterns in textured images and spectrograms and is the origin for the operations by \cite{phan2013basis, eikmeier2018hyperkron, wang2021kronecker, colley2023dominant}; Ragnarsson provides the most comprehensive discussion of the tensor KP, in the context of block tensor unfolding and factorization, although with little attribution \cite[Chapter 4]{ragnarsson2012structured}. Several proposed instances of this operation, labeled either as a \KP \cite{mohammadi2016triangular, kressner2017recompression, batselier2017constructive} or under alternative names like the ``direct product" \cite{shao2013general}, have discussed a similar operator without proper attribution to earlier works, but a comprehensive list of all such instances is currently unavailable.

Continuing the story of Zehfuss and Kronecker, along with the tensor extensions by Nicholas and others, this article makes the following contributions:
\begin{itemize}
    \item We provide a comprehensive list of algebraic, structural, and spectral properties and derive tensor decompositions, such as the Tucker, Higher-Order Singular Value, Canonical Polyadic, orthogonal, and Tensor Train decompositions in terms of the tensor KP.
    \item We introduce and characterize the structure of Kronecker hypergraphs in terms of the \KP of adjacency tensors.
    Characterizations of the stability of polynomial dynamics \textcolor{black}{defined according to} Kronecker hypergraph \textcolor{black}{structure is}  expressed in terms of the factor hypergraphs.
    \item We showcase a several order of magnitude improvement in run time when computing Canonical Polyadic and Tensor Train decompositions as well as Z-eignevectors of modest sized tensors with Kronecker structure. Furthermore, we present a concrete example demonstrating the polynomial stability result.
\end{itemize}

We discuss the structural, spectral, and decomposition properties of the tensor \KP (\cref{sec:2}), examine structural and dynamic properties of Kronecker hypergraphs (\cref{sec:3}), and provide several numerical examples where Kronecker structure is utilized to calculate eigenvalues, tensor decompositions, and polynomial stability in improved time (\cref{sec numerical examples}).

\section{Tensor Kronecker Product}\label{sec:2}
In this section we investigate various algebraic, structural, and spectral properties of the tensor KP. We utilize these properties to express and compute Tucker and Tensor Train decompositions with tensor KP.

\subsection{Tensor preliminaries}
Tensors are multidimensional arrays generalized from vectors and matrices\footnote{In precise mathematical terminology, a multidimensional array is not formally equivalent to a tensor \cite[Chapter 14]{roman2005advanced}. Nevertheless, for the sake of consistency with the nomenclature used in \cite{kolda2009tensor}, we adopt the practice of referring to multidimensional arrays as tensors.}. The order of a tensor is the number of its indices, also referred to as modes. A $k$th-order real-valued tensor is often denoted by $\tT\in \mathbb{R}^{n_1\times  \dots \times n_k}$, where $n_j$ is the dimension of the $j$th mode. A tensor $\tT$ is supersymmetric if $\tT_{j_1\dots j_k}$ is invariant under any permutation of its indices. Moreover, tensor norms are defined similarly to their matrix counterparts.

\begin{definition}[Tensor norms \cite{kolda2009tensor}]
    Given a $k$th-order tensor $\tT\in \mathbb{R}^{n_1\times  \dots \times n_k}$, the $l_1$-norm and the Frobenius norm, denoted by $\|\cdot\|_1$ and $\|\cdot\|_F$, are 
    \begin{equation*}
        \|\tT\|_1=\sum_{j_1=1}^{n_1}\dots\sum_{j_k=1}^{n_k}|\tT_{j_1\dots j_k}|\text{ and }\|\tT\|_F=\sqrt{\sum_{j_1=1}^{n_1}\dots\sum_{j_k=1}^{n_k}\tT_{j_1\dots j_k}^2},
    \end{equation*}
    respectively.
\end{definition}

In order to represent the fibers and slices of a tensor (as defined below), we employ the MATLAB colon notation. For example, $\Av_{i:}$ represents the $i$th rows of a matrix $\Av$, while $\Av_{:j}$ corresponds to the $j$th columns of $\Av$.

\begin{definition}[Fibers and slices \cite{kolda2009tensor}]
Fibers refer to vectors obtained from a tensor by keeping all indices fixed except for one. Slices, on the other hand, correspond to matrices obtained from a tensor by fixing all indices except for two.
\end{definition}

Consider a third-order tensor $\tT \in \mathbb{R}^{n_1 \times n_2 \times n_3}$. The fibers $\tT_{:jl}$, $\tT_{i:l}$, and $\tT_{ij:}$ are commonly referred to as columns, rows, and tubes. The slices $\tT_{i::}$, $\tT_{:j:}$, and $\tT_{::l}$ are commonly known as transverse, lateral, and frontal slices, respectively. Fibers and slices can be further generalized by fixing any set of indices as the number of modes increases. For instance, if $\tT\in\R^{n_1\times n_2\times n_3\times n_4},$ we denote  $\tT_{4=1}$, where $4$ indicates the fixed mode and $1$ is the selected index, as the tensor $\tT_{:::1}.$ These structures play a fundamental role in the tensor decompositions of the tensor KP. 

\subsubsection{Tensor products}
Multiplying a tensor $\tT$ by a vector, matrix, or another tensor involves multiplying each fiber of $\tT$ along a specified mode with the corresponding vector, rows/columns of the matrix, or fibers of the other tensor. 

\begin{definition}[Tensor inner product \cite{kolda2009tensor}]
Given two $k$th-order tensors $\tT\in \mathbb{R}^{n_1\times  \dots \times n_k}$ and $\tS\in \mathbb{R}^{n_1\times  \dots \times n_k}$ of the same size, the tensor inner product is  
\begin{equation*}
    \langle \tT, \tS\rangle = \sum_{j_1=1}^{n_1}\dots \sum_{j_k=1}^{n_k} \tT_{j_1\dots j_k}\tS_{j_1\dots j_k}.
\end{equation*}
\end{definition}

The tensor Frobenius norm therefore can be expressed as $\|\tT\|^2_F=\langle \tT, \tT\rangle$. Two tensors $\tT$ and $\tS$ are called orthogonal if $\langle \tT, \tS\rangle=0$. The tensor outer product is a direct extension of the vector outer product.

\begin{definition}[Tensor outer product \cite{kolda2009tensor}]
Given a $k_1$th-order tensor $\tT\in \mathbb{R}^{n_1\times  \dots \times n_{k_1}}$ and a $k_2$th-order tensor $\tS\in \mathbb{R}^{m_1\times  \dots \times m_{k_2}}$, the tensor outer product is
\begin{equation*}
    (\textsf{T}\circ\textsf{S})_{j_1\dots j_{k_1} i_1\dots i_{k_2}} = \textsf{T}_{j_1\dots j_{k_1}}\textsf{S}_{i_1\dots i_{k_2}}.
\end{equation*}
\end{definition}

Furthermore, tensors can be multiplied by vectors or matrices, similar to matrix vector or matrix matrix multiplications.

\begin{definition}[Tensor matrix multiplication \cite{kolda2009tensor}]\label{def: tmm}
Given a $k$th-order tensor $\tT\in \mathbb{R}^{n_1\times  \dots \times n_k}$, the tensor matrix multiplication $\tT \times_{p} \Av\in\R^{n_1\times \dots\times n_{p-1}\times m\times n_{p+1}\times \dots\times n_k}$ along mode $p$ for a matrix $\Av\in\R^{m\times n_p}$ is 
\begin{equation*}
(\tT \times_{p} \Av)_{j_1\dots j_{p-1}ij_{p+1}\dots j_k}=\sum_{j_p=1}^{n_p}\tT_{j_1\dots j_p\dots j_k}\Av_{ij_p}.
\end{equation*}
\end{definition}

The tensor vector multiplication follows directly from the tensor matrix multiplication by treating the vector as a matrix. Furthermore, the tensor vector multiplication can be extended to 
\begin{equation}\label{eq6}
\begin{split}
\tT\times_1 \textbf{v}_1 \times_2\dots \times_{k}\textbf{v}_k\in\R,
\end{split}
\end{equation}
for $\textbf{v}_p\in \R^{n_p}$. The above expression is also known as the homogeneous polynomial associated with $\tT$. If $\textbf{v}_p=\textbf{v}$ for all $p$, we write  (\ref{eq6})  as $\tT\textbf{v}^{k}$ for simplicity. Last but not least, the Einstein product is one particular form of tensor tensor multiplications.

\begin{definition}[Einstein product \cite{cui2016eigenvalue}]\label{def: tem}
Given an even-order tensor\\
$\tT\in\R^{n_1\times m_1\times\dots\times n_k\times m_k}$ and a $k$th-order tensor $\tX\in\R^{m_1\times\dots\times m_k},$ the Einstein product $\tT*\tX\in\mathbb{R}^{n_1\times \dots\times n_k}$ is 
\begin{equation*}
    (\tT*\tX)_{j_1\dots j_k}=\sum_{i_1=1}^{m_1}\dots\sum_{i_k=1}^{m_k}\tT_{j_1i_1\dots j_ki_k}\tX_{i_1\dots i_k}.
\end{equation*}
\end{definition}

The above tensor products are fundamental to tensor eigenvalue problems and tensor decompositions and will be used throughout the remainder of the article.

\subsubsection{Spectral definitions}
Introduced independently in 2005 by Qi \cite{qi2005eigenvalues} and Lim \cite{lim2005singular}, various versions of tensor eigenvalues have since been developed. Tensor eigenvalues extend the concept from matrices to account for repeated tensor vector or tensor tensor multiplications.

\begin{definition}[H-eigenpair \cite{qi2005eigenvalues}]
    Given a $k$th-order supersymmetric  tensor $\tT\in\mathbb{R}^{n\times \dots\times n}$, the H-eigenpair ($\lambda,\xv$), where $\lambda\in\R$ and $\xv\in\R^n,$ is a solution to \label{def: Heigen}
    \begin{equation}\label{def:eq: H-eigenpair}
        \tT\xv^{k-1}=\lambda\xv^{\{k-1\}},
    \end{equation}
    where $\xv^{\{k-1\}}$ denotes Hadamard product exponentiation to the $(k-1)$th power. 
\end{definition}

\begin{definition}[Z-eigenpair \cite{qi2005eigenvalues}]
    Given a $k$th-order supersymmetric tensor $\tT\in\mathbb{R}^{n\times \dots\times n}$, the Z-eigenpair $(\lambda,\xv),$ where $\lambda\in\R$ and $\xv\in\R^n,$  is a solution to \label{def: Zeigen}
    \begin{equation}\label{def:eq: Z-eigenpair}
        \tT\xv^{k-1}=\lambda\xv\text{ and } \xv^\top\xv=1.
    \end{equation}
\end{definition}

\begin{definition}[Generalized M-eigentriple \cite{qi2009conditions}]
Given an even-order ($2kth$-\newline order) tensor $\tT\in\mathbb{R}^{n\times n \dots\times n\times n}$, the $M$-eigentriple $(\lambda,\xv,\yv),$ where $\lambda\in\R$, $\xv$ and $\yv\in\R^n$, and $\xv^\top\xv=\yv^\top\yv=1$, is the solution to
\begin{equation*}
    \tT\xv^{k}\yv^{k-1}=\lambda\yv\text{ and }\tT\yv^{k}\xv^{k-1}=\lambda\xv.
\end{equation*}
\end{definition}

\begin{definition}[U-eigenpairs \cite{cui2016eigenvalue}]
Given an even-order ($2kth$-order) tensor \newline $\tT\in\mathbb{R}^{n_1\times n_1 \times \dots\times n_k\times n_k}$, the U-eigenpair $(\lambda, \tX),$ where  $\lambda\in\C$ and $\tX\in\C^{n_1\times \dots\times n_k}$, is the solution to \label{def: Ueigen}
\begin{equation*}
    \tT*\tX=\lambda\tX.
\end{equation*}
\end{definition}

The computation of tensor eigenvalues continues to be a topic of active research. By leveraging  mixed product properties, we can utilize the Kronecker structure of a tensor to facilitate the computation of specific tensor eigenvalues.

\subsubsection{Tensor decompositions}
A wide range of tensor decompositions have been developed \cite{kolda2009tensor, oseledets2011tensor}. These decompositions serve as powerful tools to represent and analyze complex data structures that exhibit multidimensional relationships.

\begin{definition}[Tucker decomposition \cite{tucker1966some}]
The Tucker decomposition represents a $k$th-order tensor  $\tT\in\R^{n_1\times\dots\times n_k}$ as a core tensor $\tS$ multiplied by matrices $\Uv_p$ along each mode $p$, i.e., 
\begin{equation}\label{def: eq: tucker decomposition}
    \tT = \tS\times_1\Uv_1\times_2\dots\times_k\Uv_k,
\end{equation}
where $\tS\in\R^{m_1\times\dots\times m_k}$ and $\Uv_p\in\R^{n_p\times m_p}$.
\end{definition}

The Canonical Polyadic decomposition (CPD),  orthogonal decomposition, and  higher-order SVD (HOSVD) are special cases of the Tucker decomposition. In the case of CPD, $\tS$ becomes a diagonal tensor such that  $\tT$ can be equivalently represented as the sum of rank one tensors, and the number of nonzero diagonal entries in $\tS$ is referred to as the CP-rank \cite{kiers2000towards}. For orthogonal decomposition,  a supersymmetric tensor is called orthogonally decomposable when $\tS$ is diagonal and the matrices $\Uv_j$ contain orthanormal columns \cite{robeva2016orthogonal}. The HOSVD is a multilinear extension of the matrix SVD that requires the matrices $\Uv_j$ to contain orthonormal columns and the core tensor \textsf{S} to satisfy the following two conditions \cite{lathauwer2000multilinear}:
\begin{itemize}
    \item All-orthogonality: the sub-tensors $\tS_{j_p=\alpha}$ and $\tS_{j_p=\beta}$ are orthogonal for all possible values $p,\alpha,\beta$ subject to $\alpha\neq\beta.$
    \item Ordering: $\|\tS_{j_p=1}\|_F\geq\dots\geq\|\tS_{j_p=n_p}\|_F\geq 0$ for all possible values of $p$.
\end{itemize}
The Frobenius norms $\|\tS_{j_p=\alpha}\|_F=\gamma_\alpha^{(p)}$ are the $p$-mode singular values of $\tT,$ and the number of nonvanishing $p$-mode singular values is called the $p$-rank of $\textsf{T}$. The multilinear rank of $\tT$ is defined as the maximum of its $p$-ranks.

The Tensor Train decomposition (TTD) is an alternative decomposition to Tucker decompositions.

\begin{definition}[TTD \cite{oseledets2010tt, oseledets2011tensor}]
Given a $k$th-order tensor $\tT\in\R^{n_1\times \dots\times n_k}$, the TTD represents $\tT$ as a sum of the outer product of fibers taken from third-order tensors, i.e., 
\begin{equation}
\tT=\sum_{j_0=1}^{r_0}\dots\sum_{j_k=1}^{r_k}\tT^{(1)}_{j_0:j_1}\circ\dots\circ\tT^{(k)}_{j_{k-1}:j_k},
\end{equation}
where $\{r_0,\dots,r_k\}$ is called the set of  TT-ranks and $\tT^{(p)}\in\R^{r_{p-1}\times n_p\times r_p}$ are called the core tensors of $\tT$.
\label{def: TTD}
\end{definition}

The computation of the TTD is numerically stable, exhibiting linear complexity in both time and storage relative to the dimensions of $\tT$. This efficiency makes it particularly suitable for addressing high-dimensional problems.

\subsection{Block Definition of KP}
Block tensors, developed using principles similar to block matrix operations, are intimately related to the tensor KP. Just like block matrices, a block tensor is a tensor in which the elements themselves are tensors. In the case of a $k$th-order block tensor $\tT$, we denote the $(j_1,\dots,j_k)$th block of $\tT$ as $[\tT]_{j_1\dots j_k}$ where brackets denote tensor blocks. As outlined in \cite{ragnarsson2012block}, block tensors have three principle advantages: (1) \textit{structure}: block-level sparsity is a common feature of many data and algorithms; (2) \textit{generalizability}: many operations naturally extend to block tensors; (3) \textit{performance}: block operations are often faster than elementwise methods, as we will see in \cref{sec numerical examples}.

Similar to matrices, operations on block tensors can be defined and executed either elementwise or through vectorization, bypassing the use of blocks. This viewpoint was adopted by \cite{nicholson2001kronecker1, akoglu2008rtm, phan2012revealing, batselier2017constructive} in their independent introductions of the tensor KP. Using block tensors, we can alternatively define the tensor \KP in a manner more consistent with the matrix KP.

\begin{definition}[Tensor Kronecker product]
    Given two $k$th-order tensors $\tB\in\R^{n_1\times\dots\times n_k}$ and $\tC\in\R^{m_1\times\dots\times m_k}$, the tensor Kronecker product  $\tB\otimes\tC\in\R^{n_1m_1\times\dots\times n_km_k}$ is defined as  \label{def: tensor Kronecker product}
    \begin{equation}\label{def: eq: tensor kron}
        [\tB\otimes\tC]_{i_1\dots i_k}=\tB_{i_1\dots i_k}\tC.
    \end{equation}
\end{definition}

In accordance with the dimensions of $\tB$ and $\tC$ provided in \cref{def: tensor Kronecker product}, the tensor \KP can also be defined from an element-wise perspective as 
\begin{equation}\label{eq: kron tensor index mapping inv}
    (\tB\otimes \tC)_{((i_1-1)m_1+j_1)((i_2-1)m_2+j_2)\dots((i_k-1)m_k+j_k)}=\tB_{i_1\dots i_k}\tC_{j_1\dots j_k},
\end{equation}
similar to the products used in \cite{akoglu2008rtm, phan2012revealing}; the closest definition is proposed in \cite{ragnarsson2012structured}, where linear, rather than tuple, indices are used. The adoption of block tensors in \cref{def: tensor Kronecker product} allows for a more natural development of various properties for the tensor KP, as opposed to an element-wise definition. \cref{fig: Kronecker block tensor} provides a visual representation of the tensor KP, using the concept of block tensors.

\subsection{Algebraic properties}\label{sec: tkp algebra}
The tensor \KP exhibits various algebraic properties similar to those of  the matrix KP. Just like the matrix KP, the tensor KP $\otimes$ is bilinear, meaning that the function$f(\tB,\tC)=\tB\otimes\tC$ is linear with respect to both $\tB$ and $\tC$.

\begin{property}[Kronecker product algebra 
\cite{batselier2017constructive}]
Given three $k$th-order tensors $\tB$, $\tC$, and $\tD$ with a scalar $\alpha$, the following relations hold:
\label{prop: kronecker product algebra}
\begin{itemize}
    \item Distribution over tensor addition:
    \begin{equation}
        \tB\otimes(\tC+\tD) = \tB\otimes\tC + \tB\otimes\tD \text{ and } (\tB+\tC)\otimes\tD = \tB\otimes\tD + \tC\otimes\tD;
    \end{equation}
    \item Commutation with scalar multiplication:
    \begin{equation}
        \alpha(\tB\otimes\tC) = (\alpha\tB)\otimes\tC = \tB \otimes (\alpha\tC);
    \end{equation}
    \item Kronecker product association:
    \begin{equation}
        (\tB\otimes\tC)\otimes\tD = \tB\otimes(\tC\otimes\tD).
    \end{equation}
\end{itemize}
\end{property}

For the matrix KP, the mixed product property, i.e., $(\Av\Bv)\otimes(\Cv\Dv)=(\Av\otimes\Cv)(\Bv\otimes\Dv)$ \cite[Chapter 12.3.1]{golub2013matrix}, simplifies the expression of diverse matrix decompositions. Similarly, translating tensor decompositions into KPs necessitates mixed product properties, each corresponding to the various tensor products involved.

\begin{property}[Mixed tensor matrix multiplication 
{\cite[Chapter 4]{ragnarsson2012structured}}]
    Given two $k$th-order tensors $\tB\in\R^{n_1\times\dots\times n_k}$ and $\tC\in\R^{m_1\times\dots\times m_k}$, the following relation holds:
    \begin{equation}
        (\tB\otimes\tC)\times_p(\Xv\otimes\mathbf{Y}) = (\tB\times_p\Xv)\otimes(\tC\times_p\mathbf{Y})
    \end{equation}
    for any two matrices $\Xv\in\R^{r\times n_p}$ and $\mathbf{Y}\in\R^{h\times m_p}$.
\label{prop: mixed product tensor matrix}
\end{property}

Properties \ref{prop: kronecker product algebra} and \ref{prop: mixed product tensor matrix}, are explicitly outlined in \cite{batselier2017constructive} and {\cite[Chapter 4]{ragnarsson2012structured}} respectively, and these properties are clearly derived from the index-free perspective of tensor products presented in \cite[Chapter 14, 359-378]{roman2005advanced}. The subsequent properties, namely Properties \ref{prop: mixed tensor inner product}-\ref{prop: mixed product Einstein}, follow a similar index-free perspective and are provided without proof for brevity.

\begin{property}[Mixed tensor inner product {\cite[Chapter 4]{ragnarsson2012structured}}]\label{prop: mixed tensor inner product}
    Given two $k$th-order tensors $\tB\in\R^{n_1\times\dots\times n_k}$ and $\tC\in\R^{n_1\times\dots\times n_k}$ of the same size, the following relation holds:
    \begin{equation}
        \langle \tB\otimes \tC, \tX\otimes \tY\rangle = 
        \langle \tB,\tX\rangle \langle \tC,\tY\rangle
    \end{equation}
    for any two tensors $\tX\in\R^{n_1\times\dots\times n_k}$ and $\tY\in\R^{n_1\times\dots\times n_k}$.
\label{prop: mixed inner product}
\end{property}

\begin{property}[Mixed tensor outer product]
    Given two $k_1$th-order tensors $\tB\in\R^{n_1\times\dots\times n_{k_1}}$ and  $\tC\in\R^{m_1\times\dots\times m_{k_1}}$, the following relation holds:
    \begin{equation}
        (\tB\otimes\tC)\circ(\tX\otimes\tY) = (\tB\circ\tX)\otimes(\tC\circ\tY)
    \end{equation}
    for any two $k_2$th-order tensors $\tX\in\R^{h_1\times\dots\times h_{k_2}}$ and $\tY\in\R^{r_1\times\dots\times r_{k_2}}$.
\label{prop: mixed outer product}
\end{property}

\begin{property}[Mixed Einstein product]
    Given two even-order tensors $\tB\in\R^{n_1\times m_1\times \dots\times n_k\times m_k}$ and $\tC\in\R^{h_1\times r_1\dots\times h_k\times r_k}$, the following relation holds:
    \begin{equation}
        (\tB\otimes\tC)*(\tX\otimes\tY) = (\tB*\tX)\otimes(\tC*\tY)
    \end{equation}
    for any two $k$th-order tensors $\tX\in\R^{m_1\times\dots\times m_k}$ and $\tY\in\R^{r_1\times \dots\times r_k}$.
\label{prop: mixed product Einstein}
\end{property}

Properties \ref{prop: mixed tensor inner product}, \ref{prop: mixed outer product}, and \ref{prop: mixed product Einstein} play a crucial role in understanding and computing Kronecker tensor eigenvalues and tensor decompositions. In addition, the norms of the tensor \KP can be computed from the individual tensor norms.

\begin{property}[Separable norms]\label{prop: separable norms}
    Given two $k$th-order tensors $\tB\in\R^{n_1\times\dots\times n_k}$ and $\tC\in\R^{m_1\times\dots\times m_k},$ the $l_1$-norm and Frobenius norm of the tensor \KP $\tB\otimes\tC$ can be computed as
    \begin{equation}
        \|\tB\otimes\tC\|_1=\|\tB\|_1\|\tC\|_1\text{ and }\|\tB\otimes\tC\|_F=\|\tB\|_F\|\tC\|_F,
    \end{equation}
    respectively.
\label{prop: l1 norm}
\end{property}

For brevity, the proof is omitted. Furthermore, the tensor \KP preserves fiber, slide, and sub-tensor structures, proving crucial in computing the TTD or HOSVD of a tensor generated via this product.

\begin{property}[Fiber structure]
Given two $k$th-order tensors $\tB\in\R^{n_1\times\dots\times n_k}$ and $\tC\in\R^{m_1\times\dots\times m_k},$ the \KP of the $p$-mode fibers of $\tB$ and $\tC$ is the $p$-mode fiber of the tensor \KP $\tB\otimes\tC$.
\label{prop: kronecker fibers}
\end{property}
\begin{proof}
    Consider the $p$-mode fibers  $\tB_{i_1\dots i_{p-1}:i_{p+1}\dots i_k}$ and $\tC_{j_1\dots j_{p-1}:j_{p+1}\dots j_k}$. From \cref{eq: kron tensor index mapping inv}, the \KP of these fibers can be written as
    \begin{equation*}
    \begin{split}
        (\tB\otimes\tC)_{((i_1-1)m_1+j_1)\dots((i_{p-1}-1)m_{p-1}+j_{p-1}):((i_{p+1}-1)m_{p+1}+j_{p+1})\dots((i_k-1)m_k+j_k)},
    \end{split}
    \end{equation*}
    which happens to be a $p$-mode fiber in $\tB\otimes\tC$ as desired.
\end{proof}

\begin{property}[Slice and sub-tensor structure]\label{prop: kronecker slices}
Given two $k$th-order tensors $\tB\in\R^{n_1\times \dots\times n_k}$ and $\tC\in\R^{m_1\times\dots\times m_k},$ the \KP of the $(p_1,p_2)$-mode slices of $\tB$ and $\tC$ is a $(p_1,p_2)$-mode slice of the tensor \KP $\tB\otimes\tC$. Moreover, the \KP of sub-tensors $\tB_{j_p=\alpha}$ and $\tC_{j_p=\beta}$ forms the sub-tensor $(\tB\otimes\tC)_{j_p=(\alpha-1)m_p+\beta}.$
\end{property}

Note that Properties \ref{prop: kronecker fibers} and \ref{prop: kronecker slices} are true only if the fibers (slices or sub-tensors) are taken in the same mode(s) of a tensor.

\subsection{Tensor Kronecker structure}
Many tensor structures, such as diagonal, banded, block, stochastic \cite{benson2017spacey}, are similar to their corresponding matrix structures; for instance, a tensor is diagonal if $\tT_{j_1\dots j_k}=0$ except possibly when $j_1=\dots=j_k$. These structures are preserved under the tensor KP.

\begin{property}[Tensor Kronecker structure]
The tensor \KP of $k$th-order \label{prop: Kronecker product preserves tensor structure}
\begin{equation*}
\begin{pmatrix}\text{diagonal}\\\text{banded}\\\text{supersymmetric}\\\text{upper (lower) triangular}\\\text{stochastic}
\\\text{Hankel \cite{batselier2017constructive}}\\\text{Toeplitz \cite{batselier2017constructive}}\end{pmatrix} \text{ tensors is a }\begin{pmatrix}\text{diagonal}\\\text{block banded}\\\text{supersymmetric}\\\text{upper (lower) triangular}\\\text{stochastic}
\\\text{Hankel \cite{batselier2017constructive}}\\\text{Toeplitz \cite{batselier2017constructive}}\end{pmatrix} \text{ tensor.}
\end{equation*}
\end{property}
\begin{proof}
    These properties are relatively straightforward to verify, and as an example, we illustrate the preservation of diagonal structure. Suppose $\tB$ and $\tC$ are diagonal tensor such that $\tB_{i_1,\dots,i_k}\neq 0$ and $\tC_{j_1,\dots,j_k}\neq 0$ if and only if $i_1=\dots=i_k$ and $j_1=\dots=j_k$. It follows that $[\tB\otimes\tC]_{i_1,\dots,i_k}\neq0$ if and only if $i_1=\dots=i_k.$ More over, within the block $[\tB\otimes\tC]_{i_1,\dots,i_k}$, the $(j_1,\dots,j_k)$th element is nonzero if and only if $j_1=\dots=j_k.$ Therefore, an element of $\tB\otimes\tC$ is nonzero if and only if the element is a diagonal element of a diagonal block of $\tB\otimes\tC,$ which is equivalent to saying the element is along the main diagonal of the \KP $\tB\otimes\tC$. Preservation of the remaining structures may be verified similarly.
\end{proof}

Additionally, the tensor \KP preserves the properties of being general symmetric, persymmetric, and centrosymmetric \cite{batselier2017constructive}.

\subsection{Spectral properties}
In contrast to matrix eigenvalues, the preservation of positive or negative definiteness of tensors under the \KP is not guaranteed. Still, some tensor H-, Z-, M-, and U-eigenpairs of the tensor \KP can be computed from factor tensors.

\begin{property}[H-eigenpairs]
    Given the H-eigenpairs $(\alpha,\xv)$ and $(\beta,\yv)$ of two $k$th-order supersymmetric tensors $\tB$ and $\tC$, respectively, the pair $(\alpha\beta,\xv\otimes\yv)$ is an H-eigenpair of the tensor \KP $\tB\otimes\tC$.
    \label{prop: H-eigenpair} 
\end{property}
\begin{proof}
    According to  \cref{prop: mixed product tensor matrix}, taking the tensor \KP of the left and right sides of (\ref{def:eq: H-eigenpair}) yields
    \begin{equation*}
        \begin{split}
        \tB\xv^{k-1}\otimes\tC\yv^{k-1} & =\alpha\xv^{\{k-1\}}\otimes\beta\yv^{\{k-1\}}\Longleftrightarrow
        (\tB\otimes\tC)(\xv\otimes\yv)^{k-1} = \alpha\beta(\xv\otimes\yv)^{\{k-1\}}.
    \end{split}
    \end{equation*}
    Therefore, based on  \cref{def: Heigen}, $(\alpha\beta,\xv\otimes\yv)$ is an H-eigenpair of the tensor \KP $\tB\otimes\tC$.
\end{proof}

Following steps similar to as in the proof above, one can  establish analogous properties for other tensor eigenvalue problems. 

\begin{property}[Z-eigenpairs {\cite[Chapter 4]{ragnarsson2012structured}}]
    Given the Z-eigenpairs  $(\alpha,\xv)$ and $(\beta,\yv)$ of two $k$th-order tensors $\tB$ and $\tC$, respectively, the pair $(\alpha\beta,\xv\otimes\yv)$ is an Z-eigenpair of the tensor \KP $\tB\otimes\tC$.
    \label{prop: Z-eigenpair}
\end{property}

\begin{property}[M-eigentriples]\label{prop: M-eigentriple}
Given the M-eigentriples $(\alpha,\wv,\xv)$ and \newline $(\beta,\yv,\zv)$ of two even-order tensors $\tB$ and $\tC,$ the triple $(\alpha\beta,\wv\otimes\yv,\xv\otimes\zv)$ is a $M$-eigenpair of the tensor \KP $\tB\otimes\tC$.
\end{property}

\begin{property}[U-eigenpairs]
Given the U-eigenpairs $(\alpha,\tX)$ and $(\beta,\tY)$ of two even-order tensors $\tB$ and $\tC$, respectively, the pair $(\alpha\beta,\tX\otimes\tY)$ is a U-eigenpair of the tensor \KP $\tB\otimes\tC$.
\label{prop: U-eigenpair}
\end{property}

The converses of  Properties \ref{prop: H-eigenpair}-\ref{prop: M-eigentriple} do not hold universally. Specifically, not all eigenpairs of $\tB\otimes\tC$ will correspond to KPs of eigenpairs of $\tB$ and $\tC$. This discrepancy arises as a result of the number of eigenvalues associated with a tensor \cite{cartwright2013number}. It is worth noting, that the dominant Z-eigenpair can be expressed as the \KP of the dominant Z-eigenpairs of the factor tensors \cite{colley2023dominant}.

\subsection{Tensor decompositions}\label{sec: tensor KP decompositions}
Computing tensor decompositions, particularly for high-dimensional, high-mode tensors, can be facilitated by leveraging the \KP or tensors generated with this structure. In this subsection, we explore the use of Kronecker factors in expressing different tensor decompositions.  Our focus centers around the problem of efficiently computing the tensor decomposition of the tensor \KP $\tB \otimes \tC$ where the decompositions of $\tB$ and $\tC$ are known.

\begin{property}[Tucker decomposition]
    Given two $k$th-order tensors $\tB$ and $\tC$ with Tucker decompositions that have core tensors $\tS$ and $\textsf{R}$ and factor matrices $\Uv_p$ and $\Vv_p$, respectively, the Tucker decomposition of the tensor \KP $\tB\otimes\tC$ has a core tensor $\tS\otimes\textsf{R}$ and factor matrices $\Uv_p\otimes\Vv_p$ for all $p$.
\label{prop: kron tucker}
\end{property}
\begin{proof}
    According to  \cref{prop: mixed product tensor matrix}, taking the tensor \KP of the Tucker decompositions of $\tB$ and $\tC$ yields
    \begin{equation*}
        \begin{split}
    (\tB\otimes\tC) &= (\tS\times_1\Uv_1\times_2\dots\times_k\Uv_k)\otimes(\textsf{R}\times_1\Vv_1\times_2\dots\times_k\Vv_k)\\
    & = (\tS\otimes\textsf{R})\times_1(\Uv_1\otimes\Vv_1)\times_2\dots\times_k(\Uv_k\otimes\Vv_k).
    \end{split}
    \end{equation*}
    In this form, $(\tB\otimes\tC)$ is expressed as a core tensor $(\tS\otimes\textsf{R})$ times factor matrices $(\Uv_p\otimes\Vv_p)$ for all $p.$ This constitutes a Tucker decomposition of $\tB\otimes\tC.$
\end{proof}

Since the CPD and orthogonal decomposition are special cases of the Tucker decomposition, the Properties \ref{prop: kron cpd} and \ref{prop: kronecker tensor orthognally decomposable} follow similarly, and we state them without proof for brevity. 

\begin{property}[CPD]
    Given two $k$th-order tensors $\tB$ and $\tC$ with CPDs that have diagonal core tensors $\tS$ and $\textsf{R}$ and  factor matrices $\Uv_p$ and $\Vv_p$ with unit column length, respectively, the CPD of the tensor \KP $\tB\otimes\tC$ has a core tensor $\tS\otimes\textsf{R}$ and factor matrices $\Uv_p\otimes\Vv_p$ for all $p$.
\label{prop: kron cpd}
\end{property}

The CPD requirement that the core tensor $\tS$ from the Tucker decomposition (\ref{def: eq: tucker decomposition}) is diagonal and remains diagonal when computing the CPD of $\tB\otimes\tC$ is fulfilled by \cref{prop: Kronecker product preserves tensor structure}. A similar factorization can be shown where the CPD is stated in terms of a summation over matrix outer products \cite[Chapter 4]{ragnarsson2012structured}. The diagonal requirement on the core tensor is similarly applied to the orthognal decomposition, and the additional constraint that the factor matrices of orthogonal decomposition remain orthogonal is satisfied by the well known property of the matrix \KP which preserves orthogonality.

\begin{property}[Orthogonal decomposition]\label{prop: kronecker tensor orthognally decomposable}
    Given two $k$th-order supersymmetric tensors $\tB$ and $\tC$ with orthogonal decompositions that have diagonal core tensors $\tS$ and $\textsf{R}$ and orthogonal factor matrices $\Uv$ and $\Vv$, respectively, the orthogonal decomposition of the tensor \KP $\tB\otimes\tC$, which is supersymmetric, has a core tensor $\tS\otimes\textsf{R}$ and factor matrix $\Uv\otimes\Vv$.
\label{prop: kron od}
\end{property}

Whereas the CPD and orthogonal decompositions of $\tB\otimes\tC$ come directly from the decompositions of $\tB$ and $\tC,$ when $\tC$ is not sufficiently well conditioned, computing the HOSVD of $\tB\otimes\tC$ can require reindexing the separate HOSVDs of $\tB$ and $\tC$.

\begin{property}[HOSVD]
    Given two $k$th-order tensors $\tB$ and $\tC$ with HOSVDs that have core tensors $\tS$ and $\tR$ and orthogonal factor matrices $\Uv_p$ and $\Vv_p$, respectively, the HOSVD of the tensor \KP $\tB\otimes\tC$ has a core tensor that is generated by permuting the indices of $\tS\otimes\tR$ and factor matrices $\Uv_p\otimes\Vv_p$ with similarly permuted indices for all $p$.
\label{prop: kron hosvd}
\end{property}

\begin{proof}
Consider two $k$th-order tensors $\tB$ and $\tC$ with HOSVDs that have core tensors $\tS\in\R^{n_1\times\dots\times n_k}$ and $\textsf{R}\in\R^{m_1\times\dots\times m_k}$ and orthogonal factor matrices $\Uv_p$ and $\Vv_p$, respectively. In order for $\tS\otimes\textsf{R}$ and factor matrices $\Uv_p\otimes\Vv_p$ for all $p$ to constitute a valid HOSVD, both the all-orthogonality and ordering properties must be satisfied by $\tS\otimes\tR:$
\begin{itemize}
    \item All-orthogonality: From the HOSVD of $\tB$ and $\tC,$ the all-orthogonality property stipulates $\langle\tS_{j_p=\alpha_1},\tS_{j_p=\beta_1}\rangle=0$ and $\langle\tR_{j_p=\alpha_1},\tR_{j_p=\beta_1}\rangle=0$ for all valid $p,\alpha_1,\alpha_2,\beta_1,\beta_2,$ subject to $\alpha_1\neq\beta_1$ and $\alpha_2\neq\beta_2.$ The sub-tensor of $(\tS\otimes\tR)_{j_p=\alpha_3}$ can be expressed as the product of sub-tensors $\tS_{j_p=\alpha_1}\otimes\tR_{j_p=\alpha_2},$ where $\alpha_1$ and $\alpha_2$ are set according to \cref{prop: kronecker slices}. From \cref{prop: mixed tensor inner product}, the inner product of any two sub-tensors of $\tS\otimes\tR$ can be written
    \begin{equation*}
    \begin{split}
        \langle(\tS\otimes\tR)_{j_p=\alpha_3},(\tS\otimes\tR)_{j_p=\beta_3}\rangle & =\langle \tS_{j_p=\alpha_1}\otimes\tR_{j_p=\alpha_2},\tS_{j_p=\beta_1}\otimes\tR_{j_p=\beta_2}\rangle\\
        & =\langle\tS_{j_p=\alpha_1},\tS_{j_p=\beta_1}\rangle\langle\tR_{j_p=\alpha_2},\tR_{j_p=\beta_2}\rangle = 0.
    \end{split}
    \end{equation*}
    All mode $j_p$-mode sub-tensors of $\tS\otimes\tR$ are orthogonal, so $\tS\otimes\tR$ obeys the all-orthogonality property.
    \item Ordering: From the HOSVD of $\tB$ and $\tC,$ the ordering principle requires $\|\tS_{j_p=1}\|_F\geq\dots\geq\|\tS_{j_p=n_p}\|_F$ and $\|\tR_{j_p=1}\|_F\geq\dots\geq\|\tR_{j_p=n_p}\|_F$ for all valid $p.$ The $j_p$-mode sub-tensors of $\tS\otimes\tR$ have the relation $\|(\tS\otimes\tR)_{j_p=\epsilon}\|_F=\|\tS_{j_p=\alpha}\|_F\|\tR_{j_p=\gamma}\|_F,$ which allows us to apply the mixed product property \ref{prop: separable norms} and rewrite the ordering constraint $\|(\tS\otimes\tR)_{j_p=1} \|_F\geq\dots\geq\|(\tS\otimes\tR)_{j_p=n_pm_p}\|_F$ in terms of the $p$-mode singular values of $\tS$ and $\tR$ as
    \begin{equation}\label{eq: ordering of p-mode SV}
    \begin{split}
        \|\tS_{j_p=1}\|_F&\|\tR_{j_p=1}\|_F\geq\dots\geq \|\tS_{j_p=1}\|_F\|\tR_{j_p=m_p}\|_F \geq\dots\geq \\
        &\|\tS_{j_p=n_p}\|_F \|\tR_{j_p=1}\|_F\geq\dots\geq \|\tS_{j_p=n_p}\|_F\|\tR_{j_p=m_p}\|_F.
    \end{split}
    \end{equation}
    When $\|\tS_{j_p=i+1}\|_F\|\tR_{j_p=m_p}\|_F\geq\|\tS_{j_p=i}\|_F\|\tR_{j_p=1}\|_F$ for all indices $i$, \cref{eq: ordering of p-mode SV} is satisfied. If $\tC$ is sufficiently well conditioned such that $\|\tR_{j_p=1}\|_F-\|\tR_{j_p=m_p}\|_F$ is small, this criteria is satisfied, and $\tS\otimes\tR$ obeys the ordering property. Otherwise, when the $p$-mode singular values of $\tC$ are not well conditioned, the $p$-mode indices of $\tS\otimes\tR$ may be rearranged so that the ordering property is satisfied. This rearrangement must coincide with the rearrangement of the indices in $\Uv_p\otimes\Vv_p$. Since the all-orthogonality principal is satisfied, rearranging the indices in one mode of $\tS\otimes\tR$ does not impact the singular values in any other mode nor does it cause the all-orthogonality principle to be violated. By performing this rearrangement of indices in all modes, the ordering principle can be satisfied, and the HOSVD of $\tB\otimes\tC$ can be generated from the HOSVDs of $\tB$ and $\tC.$
\end{itemize}
We conclude that the all-orthogonality property is always obeyed, and it is possible to permute the indices of $\tS\otimes\tR$ to satisfy the ordering.
\end{proof}

A similar instance of permuting indices used in the proof above is required when computing the SVD of the matrix $\Bv\otimes\Cv$ from the SVDs of $\Bv$ and $\Cv$ separately. The thesis of Ragnarsson investigated the all-orthogonal property of $\tS\otimes\tR$ in \cite[Chapter 4]{ragnarsson2012structured} in the context of the HOSVD. Regardless of if the singular values must be reordered, the $p$-mode singular values of the tensor \KP $\tB\otimes\tC$ are equal to the product of the $p$-mode singular values of $\tB$ and $\tC$. Furthermore, the $p$-mode rank of the tensor \KP is the product of the $p$-mode ranks of $\tB$ and $\tC.$

The multilinear rank of the tensor \KP $\tB\otimes\tC$ is less than or equal to the product of the multilinear ranks of $\tB$ and $\tC$; this is equal only when the $p$-mode ranks of $\tB$ and $\tC$ are maximized in the same mode for both tensors, which occurs in the case of supersymmetric tensors such as those associated with  hypergraphs.

\begin{property}[TTD]\label{prop: kron TTD}
Given two $k$th-order tensors $\tB$ and $\tC$ in the  TTD form with TT-ranks $r_p$ and $s_p$ and core tensors $\tB^{(p)}$ and $\tC^{(p)}$, respectively, the tensor \KP $\tB\otimes\tC$ can be computed in the TTD form with TT-ranks $r_p s_p$ and core tensors $\tB^{(p)}\otimes\tC^{(p)}$ for all $p$.
\end{property}
\begin{proof}
    Suppose that $\tB^{(p)}$ and $\tC^{(p)}$ are the core tensors of $\tB$ and $\tC$ with TT-ranks $\{r_0,\dots,r_k\}$ and $\{s_0,\dots,s_k\}$, respectively. Taking the tensor \KP of the TTDs of $\tB$ and $\tC$ yields
    \begin{equation*}
    \begin{split}
        \tB\otimes\tC & = \bigg(\sum_{i_0=1}^{r_0}\dots\sum_{i_k=1}^{r_k}\tB^{(1)}_{i_0:i_1}\circ\dots\circ\tB^{(k)}_{i_{k-1}:i_{k}}\bigg)\otimes\bigg(\sum_{j_0^\tC=1}^{s_0}\dots\sum_{j_k=1}^{s_k}\tC^{(1)}_{j_0:j_1}\circ\dots\circ\tC^{(k)}_{j_{k-1}:j_{k}}\bigg)\\
        &= \sum_{i_0=1}^{r_0}\dots\sum_{i_k=1}^{r_k}\sum_{j_0=1}^{s_0}\dots\sum_{j_k=1}^{s_k}(\tB^{(1)}_{i_0:i_1}\circ\dots\circ\tB^{(k)}_{i_{k-1}:i_{k}})\otimes(\tC^{(1)}_{j_0:j_1}\circ\dots\circ\tC^{(k)}_{j_{k-1}:j_{k}})\\
        & = \sum_{i_0=1}^{r_0}\sum_{j_0=1}^{s_0}\dots\sum_{i_k=1}^{r_k}\sum_{j_k=1}^{s_k}(\tB^{(1)}_{i_0:i_1}\otimes\tC^{(1)}_{j_0:j_1})\circ\dots\circ(\tB^{(k)}_{i_{k-1}:i_{k}}\otimes\tC^{(k)}_{j_{k-1}:j_k}).
        \end{split}
    \end{equation*}
    The third line follows as an application of \cref{prop: mixed outer product}. Based on  \cref{prop: kronecker fibers}, the above expression can be rewritten as
    \begin{equation*}
        \tB\otimes \tC = \sum_{l_0=1}^{r_0s_0}\dots\sum_{l_k=1}^{r_ks_k} (\tB^{(1)}\otimes \tC^{(1)})_{l_0:l_1}\circ \dots \circ (\tB^{(k)}\otimes \tC^{(k)})_{l_{k-1}:l_k}.
    \end{equation*}
    In this form, $(\tB^{(i)}\otimes\tC^{(i)})$ are the core tensors of a TT-decomposition of $\tB\otimes\tC,$ and from the summations, $r_is_i$ are the TT-ranks of $\tB\otimes\tC.$
\end{proof}

The expression of the Tucker Decompositions and the TTD for $\tB\otimes\tC$ in terms of the respective decompositions of $\tB$ and $\tC$ is of interest when computing the factorizations. For large systems where the Kronecker factorization is known, as we will see in \cref{sec numerical examples}, it is beneficial to compute decompositions in terms of the Kronecker factors.

\section{Kronecker hypergraphs}\label{sec:3}
Exploring the structures and dynamics of graph and hypergraph products has been a long-standing and active objective in graph theory and network science \cite{frucht1949groups, sabidussi1959composition, weichsel1962kronecker, vizing1963cartesian, leskovec2010kronecker, eikmeier2018hyperkron, colley2023dominant}, and there is growing work on tensor-free hypergraph products  \cite{hellmuth2012survey, ostermeier2012cartesian, kaveh2015hypergraph,ostermeier2015relaxed,hellmuth2016fast}.
\textcolor{black}{Recent focus directed towards exploring higher order structures, such as motifs and multilayer networks, can be seen through the lens of graph and hypergraph products \cite{eikmeier2018hyperkron, sayama2018graph}.}
In this section, we outline preliminaries on hypergraphs, introduce Kronecker hypergraphs through the tensor KP, and  explore the structural and dynamic properties of Kronecker hypergraphs.

\subsection{Hypergraph preliminaries}
A hypergraph $\h$ is a pair $\{\V,\e\}$ where $\V$ is the node set and $\e\subseteq \mathcal{P}(\V)\setminus \{\emptyset\}$ is the hyperedge set, where $\mathcal{P}(\V)$ is the power set of $\V$. 
We exclusively consider $k$-uniform, unweighted hypergraphs, where all hyperedges contain exactly $k$ nodes, throughout the remainder of the paper.

\subsubsection{Hypergraph structure}
As the multiway analog of an adjacency matrix, an adjacency tensor is the principle representation of a hypergraph.

\begin{definition}[Adjacency tensors \cite{surana2022hypergraph}]
Given a hypergraph $\h= \{\V,\e\}$ with $n$ nodes, the adjacency tensor $\tA(\h)\in\mathbb{R}^{n\times \dots\times n}$ of $\h$, which is a $k$th order $n$-dimensional supersymmetric tensor, is defined as 
\begin{equation*}
\tA(\h)_{j_1j_2\dots j_k} = \begin{cases} \frac{1}{(k-1)!} &\text{if $\{v_{j_1},v_{j_2},\dots,v_{j_k}\}\in \e$}\\ \\0&\text{otherwise}\end{cases}.
\end{equation*}
\end{definition}

The degree $d(v)$ of a node $v\in\V$ is $d(v)=|\{e\in \e\text{ s.t. }v\in e\}|$. Similar to dyadic graphs, the degree of node $v_j$ of a 
hypergraph can be computed as
\begin{equation}\label{eq:r6}
d(v_j)=\sum_{j_2=1}^n\sum_{j_3=1}^n\dots \sum_{j_k=1}^n\tA(\h)_{jj_2j_3\dots j_k},
\end{equation}
which is equivalent to the vectorized equation $\dv = \tA(\h)\mathbf{1}_n^{k-1},$ where $\mathbf{1}_n$ denotes the vector of all ones in $\R^n$ and $\dv$ is the degree vector. If all nodes have the same degree $d$, then $\h$ is $d$-regular. The degree vector, which encodes the degree distribution, is one mechanism through which the importance or centrality of nodes may be ranked.

\begin{definition}[Eigenvector centrality \cite{benson2019three}]\label{def: centrality HG}
    Given a 
    hypergraph $\h,$ an H- or Z-eigenvector centrality vector $\cv$ is any positive real vector satisfying  (\ref{def:eq: H-eigenpair}) or (\ref{def:eq: Z-eigenpair}), respectively, for the hypergraph adjacency tensor $\tA(\h)$.
\end{definition}

There are additional notions of hypergraph centrality, but we restrict ourselves to the multilinear definition of centrality for the purpose of investigating the tensor KP; see \cite{galuppi2023spectral} for more on the spectra of weighted hypergraphs. There are also several graph and matrix representations of hypergraphs.

\begin{definition}[Clique expansion]
    Given a 
    hypergraph $\h=\{\V,\e\}$ with $n$ nodes, the clique graph $\g=\{\V,\e'\}$ where $\e'=\{(v_i,v_j)\text{ s.t. }v_i,v_j\subseteq e\in\e\}.$
\end{definition}

For a hypergraph $\h,$ the adjacency matrix of the clique expansion can be written in terms of the hyperedge set $\e$ as $\Av(\g)_{ij}=|e\in\e\text{ s.t. } v_i,v_j\in e|.$ While the graph representation may reduce the computational resources required to analyze hypergraphs, the clique expansion is lossy and limits the ability of hypergraphs to unambiguously describe multiway interactions. Other graph representations, such as the star expansion, are lossless but utilize a modified vertex set from $\h.$

\subsubsection{Hypergraph dynamics}
Hypergraph dynamics can be conceptualized in various ways.
Analogous to linear dynamical systems defined according to graph structure \cite{lin1974structural}, we refer to the homogeneous polynomial dynamical system defined according to hypergraph structure as ``hypergraph dynamics".
The state vector $x(t)\in\R^n$ of a hypergraph on $n$ nodes represents the state of each node at time point $t.$

\begin{definition}[Hypergraph dynamics \cite{chen2021controllability, pickard2023observability}]\label{def: hg dynamics}
Given a 
hypergraph $\h$ with $n$ nodes, the continuous- and discrete-time dynamics of $\h$ are defined as
\begin{equation}\label{def:eq:HG dynamics C}
    \dot{\xv}(t) = \tA(\h)\xv(t)^{k-1} \text{ and }
    \xv(t+1) = \tA(\h)\xv(t)^{k-1},
\end{equation}
respectively, where $\xv(t)\in\R^n$ is the state vector. 
\end{definition}

Recent work has investigated the controllability \cite{chen2021controllability} and observability \cite{pickard2023observability, pickard2024geometric} properties of the hypergraph dynamical systems through polynomial systems theory and tensor algebra. Furthermore, the stability theory of the hypergraph dynamical systems (\ref{def:eq:HG dynamics C}) with orthogonally decomposable adjacency tensors have been developed in continuous and discrete time.

\begin{property}[Stability of continuous-time hypergraph dynamics \cite{chen2022explicit}]
Given a 
hypergraph $\h$ with $n$ nodes such that its adjacency tensor $\tA(\h)$ is orthogonally decomposable, i.e., $\tA(\h)= \sum_{j=1}^n \lambda_j \uv_j \circ\dots\circ \uv_j$, and the initial condition $\xv(0)=\sum_{j=1}^n\alpha_j\uv_j$, the equilibrium point $\xv_e=\mathbf{0}$ of the continuous-time case in \cref{def:eq:HG dynamics C} is:
\begin{itemize}
    \item stable if and only if $\lambda_j\alpha_j^{k-2}\leq 0$ for all $j=1,\dots, n$;
    \item asymptotically stable if and only if $\lambda_j\alpha_j^{k-2}< 0$ for all $j=1,\dots, n$;
    \item unstable if $\lambda_j\alpha_j^{k-2}> 0$ for some $j=1,\dots, n$.
\end{itemize}
Moreover, all other equilibrium point will inherit the behaviour of $\xv_e.$
\label{property: C stability}
\end{property}

\begin{property}[Stability of discrete-time hypergraph dynamics \cite{chen2021stability}]
Given a 
hypergraph $\h$ with $n$ nodes such that its adjacency tensor $\tA(\h)$ is orthogonally decomposable, i.e., $\tA(\h)= \sum_{j=1}^n \lambda_j \uv_j \circ\dots\circ \uv_j$, and the initial condition $\xv(0)=\sum_{j=1}^n\alpha_j\uv_j$, the equilibrium point $\xv_e=\mathbf{0}$ of the discrete-time case in \cref{def:eq:HG dynamics C} is:
\begin{itemize}
    \item stable if and only if $|\alpha_j\lambda_j^{\frac{1}{k-2}}|\leq 1$ for all $j=1,\dots,n$;
    \item asymptotically stable if and only if $|\alpha_j\lambda_j^{\frac{1}{k-2}}|< 1$ for all $j=1,\dots,n$;
    \item unstable if and only if $|\alpha_j\lambda_i^{\frac{1}{k-2}}|> 1$ for some $j=1,\dots,n$.
\end{itemize}
Moreover, all other equilibrium point will inherit the behaviour of $\xv_e.$
\label{property: D stability}
\end{property}

A $2$-uniform hypergraph is equivalent to a graph, in which the above stability properties are equivalent to the famous linear stability conditions.

\subsection{Kronecker hypergraph structure}
We define Kronecker hypergraphs in terms of the \KP of adjacency tensors and explore the structural and dynamic properties of this hypergraph product.

\begin{definition}[Kronecker hypergraphs]
    Given 
    hypergraphs $\mathcal{H}_1$ and $\mathcal{H}_2$, the Kronecker hypergraph $\mathcal{H}=\mathcal{H}_1\otimes\mathcal{H}_2$ is defined by the adjacency tensor that is the tensor \KP of the adjacency tensors of $\h_1$ and $\h_2$, \label{def: Kronecker hypergraph}
    \begin{equation}\label{def:eq: kron HG}
        \tA(\mathcal{H})=\tA(\mathcal{H}_1)\otimes\tA(\mathcal{H}_2).
    \end{equation}
\end{definition}

Note that for 2-uniform hypergraphs, which are graphs, this definition reduces to classic Kronecker graph. Since the tensor \KP does not change the number of modes for a tensor, the \KP of $k$-uniform hypergraphs remains $k$-uniform. Kronecker hypergraphs can also be defined directly in terms of their vertex and edge sets.

\begin{property}[Node and hyperedge sets]\label{prop: strong hypergraph product}
Given 
hypergraphs $\mathcal{H}_1=\{\mathcal{V}_1,\mathcal{E}_1\}$ and $\mathcal{H}_2=\{\mathcal{V}_2,\mathcal{E}_2\}$ with $m$ and $n$ hyperedges, respectively, the node and hyperedge sets of the Kronecker hypergraph $\mathcal{H}=\mathcal{H}_1\otimes\mathcal{H}_2$ are given by $\mathcal{V}=\mathcal{V}_1\times \mathcal{V}_2$ (Cartesian product) and 
\begin{equation*}
\begin{split}
    \mathcal{E}= & \big\{\{(v^{(1)}_{i_1},v^{(2)}_{j_1}),\dots,(v^{(1)}_{i_k},v^{(2)}_{j_k})\}\subseteq\mathcal{V}_1\times \mathcal{V}_2|\\
    & \{v^{(1)}_{i_1},\dots,v^{(1)}_{i_k}\}\in\mathcal{E}_1\text{ and } \{v^{(2)}_{j_1},\dots,v^{(2)}_{j_k}\}\in\mathcal{E}_2\}\big\},\\
\end{split}
\end{equation*}
respectively. 
\label{prop: V E of kron}
\end{property}
\begin{proof}
The result follows from the close correspondence between Definitions 1 and $1'$ in \cite{weichsel1962kronecker} or from Observation 1 in \cite{leskovec2010kronecker}.
\end{proof}

Based on \cref{prop: V E of kron}, the hyperedge set of $\h$ is a function of the number of hyperedges in $\e_1$,  the number of hyperedges in $\e_2$, and the number of ways a hyperedge in $\e_1$ and $\e_2$ can be aligned, of which there are $k!$ permutations. Therefore, the total number of hyperedges in the Kronecker hypergraph is $mnk!$.

Adopting the perspective of the original work on Kronecker graphs \cite{weichsel1962kronecker}, \cref{prop: strong hypergraph product} can be seen as an equivalent definition for Kronecker hypergraphs; this viewpoint has been utilized in the tensor-free graph products of Hellmuth et. al. \cite{hellmuth2012survey}. A correspondence between Kronecker hypergraphs from this vantage point.

\begin{property}[Isomorphism]\label{prop: kronecker hypergraph isomorphism}
Given 
hypergraphs $\h_1$ and $\h_2$, the Kronecker hypergraphs $\h_1\otimes\h_2$ and $\h_2\otimes\h_1$ are isomorphic.
\label{prop: kron HG isomorphism}
\end{property}
\begin{proof}
    Based on  \cref{prop: V E of kron}, there is a natural isomorphism between the node sets of $\h_1\otimes\h_2$ and $\h_2\otimes\h_1$ that preserves the hyperedge set structure.
\end{proof}

Furthermore, tensor \KP algebra from section \ref{sec: tkp algebra} reveals the degree vector, eigenvector centrality, and the clique expansion of Kronecker hyeprgraphs can be directly expressed and computed from the constituent graphs.

\begin{property}[Degree vector]
Given 
hypergraphs $\h_1$ and $\h_2$ with $m$ and $n$ nodes, respectively, with degree vectors $\dv_1$ and $\dv_2,$ the Kronecker hypergraph $\h=\h_1\otimes\h_2$ has a degree vector $\dv_1\otimes\dv_2$. Moreover, if $\h_1$ and $\h_2$ are $r$- and $s$-regular, respectively, then $\h$ is $rs$-regular.
\label{prop: kron degree sequence}
\end{property}
\begin{proof}
Applying \cref{prop: mixed product tensor matrix}, the degree vector of $\mathcal{H}$ can be computed as
\begin{equation*}
    \dv =(\tA(\h_1)\otimes\tA(\h_2))\mathbf{1}_{mn}^{k-1} = \tA(\h_1)\mathbf{1}_m^{k-1}\otimes\tA(\h_2)\mathbf{1}_n^{k-1}= \dv_1\otimes \dv_2,
\end{equation*}
where $\mathbf{1}_n$ denotes the vector of all ones in $\R^n.$ Furthermore, if $\h_1$ and $\h_2$ are $r$- and $s$-regular, i.e., all the elements in $\dv_1$ and $\dv_2$ are $r$ and $s$, respectively, then all the elements in $\dv_1\otimes\dv_2$ are $rs$, so $\h$ is $rs$-regular.
\end{proof}    

\begin{property}[Eigenvector centrality]
Given 
hypergraphs $\mathcal{H}_1$ and $\mathcal{H}_2$ with H- or Z-eigenvector centrality vectors $\cv_1$ and $\cv_2$, the H- or Z-eigenvector centrality of the Kronecker hyeprgraph $\mathcal{H}=\mathcal{H}_1\otimes\mathcal{H}_2$ is given by $\cv_1\otimes\cv_2$.
\label{prop: HKron centrality}
\end{property}
\begin{proof}
    The result follows from  the Properties \ref{prop: H-eigenpair} and \ref{prop: Z-eigenpair} in conjunction with \cref{def: centrality HG}.
\end{proof}

\begin{property}[Clique expansion]
    Given hypergraphs $\h_1$ and $\h_2$ on $n$ and $m$ vertices, respectively, the clique expansion of the Kronecker hypergraph $\h=\h_1\otimes\h_2$ is the \KP of the clique expansions of $\h_1$ and $\h_2.$
\end{property}

\begin{proof}
    The adjacency matrix of the clique expansion of $\h$ is computed as $\Av(\h)=\tA(\h)\mathbf{1}_{mn}^{k-1}.$ Applying \cref{prop: mixed product tensor matrix},
    \begin{equation*}
        \Av(\h)=\tA(\h)\mathbf{1}_{nm}^{k-2}=\tA(\h_1)\mathbf{1}_n^{k-2}\otimes\tA(\h_2)\mathbf{1}_m^{k-2}=\Av(\h_1)\otimes\Av(\h_2).
    \end{equation*}
    Given the one-to-one correspondence between adjacency tensors and hypergraphs, $\Av(\h)=\Av(\h_1)\otimes\Av(\h_2)$ indicates the clique expansion of $\h$ is the \KP of the clique expansions of $\h_1$ and $\h_2.$
    \label{prop: higher order adjacency matrices}
\end{proof}

While short, the proofs for Properties \ref{prop: kron degree sequence}, \ref{prop: HKron centrality}, and \ref{prop: higher order adjacency matrices} highlight the utility of \cref{prop: mixed product tensor matrix}. Various additional structures, such as the degree tensor, Laplacian tensor \cite{bavsic2022another}, and stochastic tensors \cite{leskovec2010kronecker, eikmeier2018hyperkron} are derived and expressed similar to their analogous structures for Kronecker graphs; however, some hypergraph operations, such as alternative graph representations do not have simple forms. For instance, the star expansion of a Kronecker hypergraph cannot be represented in terms of the star expansions of the factor hypergraphs alone.

\subsection{Kronecker hypergraph dynamics}
The dynamics of a Kronecker hypergraph can be decoupled and studied through the separate dynamics of its factor hypergraphs. For instance, given satisfactory initial conditions, the trajectory of the discrete-time equations of a Kronecker hypergraph can be integrated based upon the trajectory of the factor hypergraphs.

\begin{property}[Discrete-time trajectory]
Given hypergraphs $\mathcal{H}_1$ and $\mathcal{H}_2$ with adjacency tensors $\tA(\h)_1\in\R^{n\times\dots\times n}$ and $\tA(\h)_2\in\R^{m\times\dots\times m}$ and state vectors $\xv_1(t)\in\R^n$ and $\xv_2(t)\in\R^m$, respectively, the trajectory of the discrete-time dynamics of the Kronecker hypergraph $\mathcal{H}=\mathcal{H}_1\otimes\mathcal{H}_2$ can be computed as $\xv(t)=\xv_1(t)\otimes\xv_2(t)$ with initial condition $\xv(0)=\xv_1(0)\otimes\xv_2(0)$.
\label{thm: kron HG traj}
\end{property}

\cref{thm: kron HG traj} is verified through the application of  \cref{prop: mixed product tensor matrix}, similar to the corresponding result on Kronecker graphs \cite{chapman2014kronecker}.
The followings properties present result for the stability of Kronecker hypergraphs. 

\begin{property}[Stability of continuous-time Kronecker hypergraph dynamics]\label{prop: C kronHG stability}
Given 
hypergraphs $\h_1$ and $\h_2$ with $n$ and $m$ nodes, respectively, such that their adjacency tensors $\tA(\h_1)$ and $\tA(\h_2)$ are orthogonally decomposable, if the origin is asymptotically stable for the continuous-time factor hypergraph dynamics with initial condition $\xv_1(0)$ and $\xv_2(0)$, the origin is unstable for the continuous-time dynamics of the Kronecker hypergraph $\h=\h_1\otimes \h_2$ with the initial condition $\xv(0)=\xv_1(0)\otimes\xv_2(0)$. 
\label{prop: C kronHG stability 2}
\end{property}
\begin{proof}
    If $\tA(\h_1)$ and $\tA(\h_2)$ are orthogonally decomposable, then $\h(\tA)$ is orthogonally decomposable based on \cref{prop: kron od}. Suppose that $(\lambda_i, \vv_i)$ and $(\mu_j, \uv_j)$ are the Z-eigenpairs in the orthogonal decompositions of $\tA(\h_1)$ and $\tA(\h_2)$. Then, the initial conditions of $\h_1$ and $\h_2$ may be decomposed in terms of the orthonormal bases as $\xv_1(0)=\sum_{i=1}^n \alpha_i\vv_i$ and $\xv_2(0)=\sum_{j=1}^m \beta_j\uv_j$, respectively. The origin of $\xv(0)$ is represented as
    \begin{equation*}
        \xv(0)=\xv_1(0)\otimes\xv_2(0)=\sum_{i=1}^n\sum_{j=1}^m \alpha_i\beta_j\vv_i\otimes \uv_j,
    \end{equation*}
    where $\alpha_i\beta_j$ are coefficients of the orthonormal basis $\vv_i\otimes\uv_j.$ Since the factor hypergraph dynamical systems are  asymptotically stable at the origin, it follows that
    $\lambda_i\alpha_i^{k-2}< 0$ and $\mu_j\beta_j^{k-2}< 0$ according to \cref{property: C stability}, and from \cref{prop: Z-eigenpair}, $\lambda_i\mu_j$ are Z-eigenvalues of $\h(\tA)$ for all valid $i$ and $j.$ 
    Then,
    \begin{equation*}
       (\lambda_i\mu_j)(\alpha_i\beta_j)^{k-2}=(\lambda_i\alpha_i^{k-2})(\mu_j\beta_j^{k-2})> 0,
    \end{equation*}
    so it follows from \cref{property: C stability} that the origin is unstable on the vector field defined by $\h(\tA)$ with the initial condition $\xv(0).$
\end{proof}

\begin{property}[Stability of discrete-time Kronecker hypergraph dynamics]\label{prop: discrete-time dynamic stability}\\ 
Given 
hypergraphs $\h_1$ and $\h_2$ with $n$ and $m$ nodes, respectively, such that their adjacency tensors $\tA(\h_1)$ and $\tA(\h_2)$ are orthogonally decomposable, if the origin is asymptotically stable for the discrete-time factor hypergraph dynamics with initial condition $\xv_1(0)$ and $\xv_2(0)$, the origin is also asymptotically stable for the discrete-time dynamics of the Kronecker hypergraph $\h=\h_1\otimes \h_2$ with the initial condition $\xv(0)=\xv_1(0)\otimes\xv_2(0)$. 
\label{prop: D kronHG stability 2}
\end{property}

A proof for \cref{prop: discrete-time dynamic stability} follows a similar structure as the one of \cref{prop: C kronHG stability} but is based on \cref{property: D stability} rather than \cref{property: C stability}. Despite their similarities, Properties \ref{prop: discrete-time dynamic stability} and \ref{prop: C kronHG stability} show how the stability of Kronecker hypergraphs is not the same for discrete-time versus continuous-time systems. In discrete-time, Kronecker hypergraphs inherit both trajectory and stability properties from their factor systems. However, the continuous-time dynamics display opposite stability characteristics of the factor hypergraphs and the trajectory that is not purely separable into the trajectories of the factor hypergraphs. This is consistent with the trajectories and stability exhibited by KPs of linear dynamical systems.

\section{Numerical examples}\label{sec numerical examples}
Three examples illustrate applications of the tensor and hypergraph KPs. Subsections \ref{subsec TD} and \ref{subsec: TEC} consider the time to compute tensor decompositions and Z-eigenvalues of $\tB\otimes\tC=\tA$ with a Direct approach on $\tA$ and a Kronecker approach by performing operations on $\tB$ and $\tC$ separately. \Cref{subsec: stability example} provides an instance of continuous hypergraph dynamics to illustrate the stability of Kronecker hypergraphs (\cref{property: C stability}). All examples were performed with 16 GB of RAM, a 2.60 GHz Intel Core i7 processor, in MATLAB 2022b.\footnote{Codes utilized in \cref{sec numerical examples} is available in the software outlined in \cite{pickard2023hat} and can also be accessed here: \url{https://github.com/Jpickard1/kronecker-products-tensors-and-hypergraphs}
}

\subsection{Tensor decompositions}\label{subsec TD}
In this example, we illustrate our ability to compute the TTD and the CPD for tensors with Kronecker structure. $n$-dimensional 3rd order tensors $\tB$ and $\tC$ are randomly sampled from a uniform distribution and combined to form $\tB\otimes\tC=\tA,$ with $\tA$ containing a maximum of $2.4\times10^8$ elements.

Computing the TTD of $\tA$ requires $\mathcal{O}(kn^2r^3),$ where $k$ is the order, $n$ is the dimension of $\tB$ and $\tC$, and $r$ is an approximate TT-rank \cite{oseledets2011tensor}. The time complexity to compute the TTD can be reduced to $\mathcal{O}(kn^2r^2+kr^4)$ for a tensor with a known Tucker decomposition, and it can be further reduced for a tensor with a known Kronecker factorization. Based on \cref{prop: kron TTD}, when $\tB$ and $\tC$ are both $k$-mode $n$-dimensional tensors with approximately the same TT-ranks $r'$, computing the TTD of $\tA$ can occur in $\mathcal{O}(2knr'^3)$ where $r'\approx \sqrt{r}$. Despite relying on strong assumptions, such as the equal approximate TT-ranks of each factor tensor, the theoretical time complexity for computing the TTD of a tensor with and without a known Kronecker structure has a strong numerical agreement with our numerical experiment, as shown in \cref{fig: tensor calcs}.

The CPD is computed using a series of optimization techniques. We utilize the alternative least squares method, implemented in \cite{kolda2006matlab}, to compute the CPD of the tensor $\tA$ directly and based upon the Kronecker factors $\tB$ and $\tC$. Due to the absence of a theoretical complexity analysis of CPD computation, we present numerical comparisons to assess the runtime performance. With a fixed CP-rank, in  \cref{fig: tensor calcs}, we observe that the Kronecker approach, which computes the CPD of $\tA$ through the CPDs of $\tB$ and $\tC$, significantly outperforms the runtime of the Direct approach for various tensor sizes.

\subsection{Tensor eigenvalues}\label{subsec: TEC}
\begin{figure}[t]
    \centering
    \includegraphics[width=\textwidth]{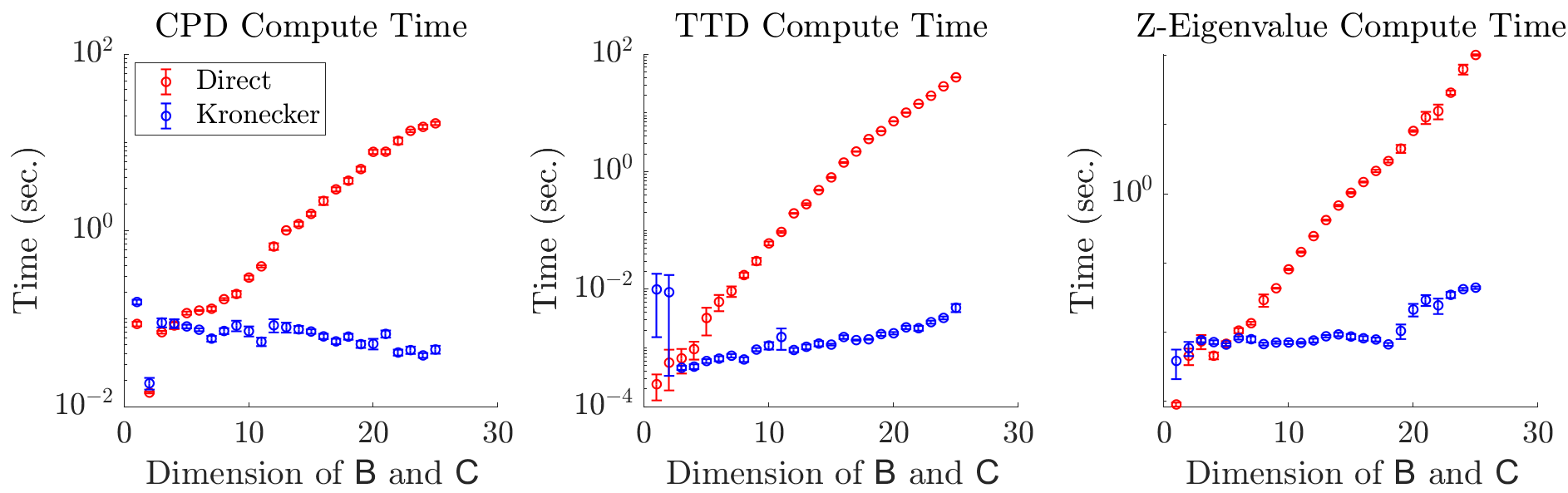}
    \vspace{-7mm}
    \caption{Runtime measurements for CPD (left), TTD (middle), and Z-eigenvalues (right) of tensors using the Direct and Kronecker approaches. The mean runtime for each tensor size is displayed over 5 trials, along with standard error bars. These calculations were performed with a CP-rank of 4, and similar results were obtained at different CP-ranks.}
    \label{fig: tensor calcs}
\end{figure}

In this example, we consider tensor eigenvalue calculations for tensors with Kronecker structure. We randomly generate two 3rd order $n$-dimensional tensors $\tB$ and $\tC$ to construct $\tB\otimes\tC=\tA$ according to the method in \cref{subsec TD}. The iterative SS-HOMP algorithm is employed to perform Z-eigenvalue calculations on the tensors $\tA,\tB,$ and $\tC$ \cite{kolda2011shifted}. Then, applying \cref{prop: Z-eigenpair} to the Z-eigenpairs of $\tB$ and $\tC,$ we obtain an alternative method of computing the Z-eigenpair of $\tA$. \cref{fig: tensor calcs} illustrates the runtime comparison for each method of calculating the Z-eigenpair of $\tA.$ For factors $\tB$ and $\tC$ with dimensions larger than 5, or when $\tA$ has a dimension larger than 25, Z-eigenvalue calculations uniformly occur faster when the Kronecker structure of $\tA$ is exploited.

The previous examples of TTD, CPD, and Z-eigenvalue maximized the numerical advantage linked to the Kronecker structure of tensor $\tA,$ where its dimensions are the square of its factors $\tB$ and $\tC.$ However, in principle, tensor decompositions, eigenvalue problems, multiplication, and other operations can be performed on tensors where the Kronecker structure does not necessarily provide evenly sized factors.

\subsection{Stability}\label{subsec: stability example}
In this example, we illustrate the counterintuitive result of  \cref{prop: C kronHG stability}, consider the two-dimensional polynomial system of degree three
\begin{equation}\label{ex: sys: stable tensor system}
    \begin{cases}
    \dot{\xv}_1 & = -1.2593\xv_1^3+1.6630\xv_1^2\xv_2-1.5554\xv_1\xv_2^2-0.1386\xv_2^3\\
    \dot{\xv}_2 & = 0.5543\xv_1^3-1.5554\xv_1^2\xv_2-0.4158\xv_1\xv_2^2-0.7036\xv_2^3
    \end{cases},
\end{equation}
borrowed from \cite{chen2021stability}. The system (\ref{ex: sys: stable tensor system}) can be represented in the form of equation \ref{def:eq:HG dynamics C} where the state transition tensor $\tB\in\R^{2\times2\times2\times2}$ is
\begin{equation}
\begin{split}
    \tB_{::11} = \begin{bmatrix}
        -1.2593&0.5534\\
        0.5543&-0.5185
    \end{bmatrix}&\text{   }
    \tB_{::12} = \begin{bmatrix}
        0.5543&-0.5185\\
       -0.5185&-0.1386
    \end{bmatrix}\\
    \tB_{::21} = \begin{bmatrix}
        0.5543&-0.5185\\
       -0.5185&-0.1386
    \end{bmatrix}&\text{   }
    \tB_{::11} = \begin{bmatrix}
        -0.5185&-0.1386\\
        -0.1386&-0.7037
    \end{bmatrix},
\end{split}
\end{equation}
and $\xv=\begin{bmatrix}
    \xv_1 & \xv_2
\end{bmatrix}^\top.$ The tensor $\tB$ is orthogonally decomposable and has strictly negative Z-eigenvalues, so the system (\ref{ex: sys: stable tensor system}) is asymptotically stable (\cref{property: C stability}). The \KP system has a state transition tensor $\tA=\tB\otimes\tB$. From \cref{prop: kron od}, $\tA$ is orthogonally decomposable, and from \cref{prop: Z-eigenpair}, $\tA$ has positive eigenvalues, so the \KP system is unstable according to \cref{prop: C kronHG stability}.\footnote{Interestingly, the tensor $\tB\otimes\tB\otimes\tB$ remains orthogonally decomposable and has all negative Z-eigenvalues, so system $\dot{\xv}=(\tB\otimes\tB\otimes\tB)\xv^3$ is stable.}

\section{Conclusion}
\textcolor{black}{
This article provides a comprehensive description of the  tensor \KP and its properties and proposes the concept of Kronecker hypergraphs as a tensor-based hypergraph product. In particular, we have demonstrated that:}
\begin{itemize}
    \item \textcolor{black}{the \KP generalizes naturally to tensors and has many nice structural, algebraic, and spectral properties in addition to being a convenient way to represent and calculate various tensor decompositions.}
    \item \textcolor{black}{hypergraph products may be expressed and understood in terms of the \KP of adjacency tensors, and the stability of homogeneous polynomial dynamics defined according to Kronecker hypergraphs may be characterized using its factored hypergraphs.}
    \item \textcolor{black}{Kronecker structure can aid large tensor decompositions and eigenvalue calculation.}
\end{itemize}
\textcolor{black}{In the future we plan to apply the tensor and hypergraph \KP concepts to real-world applications.}



\section*{Acknowledgments}{We would like to thank the two reviewers for rigorous and insightful comments. JP would also like to thank the members of the Rajapakse Lab as well as Nir Gadish and Daniel Pickard for helpful and inspiring discussion. We would also like to thank Frederick Leve at AFOSR for support and encouragement.
}

\bibliographystyle{siamplain}
\bibliography{reference}

\end{document}